\documentclass[11pt]{article}
\usepackage{amsmath, graphicx, amsfonts,amssymb, calrsfs}
\usepackage{amsfonts,mathrsfs, color, amsthm}
\allowdisplaybreaks
\usepackage{enumitem}

\def\sphere{S^{n-1}}

\def\Rn{{\mathbb R^n}}
 \def\R{\mathbb{R}}
\def\cH{\mathcal{H}}

\newtheorem{theorem}{Theorem}[section]
\newtheorem{lemma}{Lemma}[section]
\newtheorem{remark}{Remark}[section]
\newtheorem{proposition}{Proposition}[section]
\newtheorem{corollary}{Corollary}[section]
\newtheorem{prob}[theorem]{Problem}

\newtheorem{definition}{Definition}[section]

\def\bpf{\begin{proof}}
\def\epf{\end{proof}}
\def\be{\begin{equation}}
\def\ee{\end{equation}}
\def\bea{\begin{eqnarray}}
\def\eea{\end{eqnarray}}
\def\bt{\begin{theorem}}
\def\et{\end{theorem}}
\def\bl{\begin{lemma}}
\def\el{\end{lemma}}
\def\br{\begin{remark}}
\def\er{\end{remark}}
\def\bc{\begin{corollary}}
\def\ec{\end{corollary}}
\def\bd{\begin{definition}}
\def\ed{\end{definition}}
\def\bp{\begin{proposition}}
\def\ep{\end{proposition}}

\begin{document}
\title{The Minkowski problem for the non-compact convex set with an asymptotic boundary condition\footnote{The author has been supported by NSF of China (No.
		11901217 and No. 11971005)}}

\author{Ning Zhang}
\date{}
\maketitle
\begin{abstract}  

In this paper, combining the covolume, we study the Minkowski theory for the non-compact convex set with an asymptotic boundary condition. In particular, the mixed covolume of two non-compact convex sets is introduced and its geometric interpretation is obtained by the Hadamard variational formula. The Brunn-Minkowski and Minkowski inequalities for covolume are established, and the equivalence of these two inequalities are discussed as well. The Minkowski problem for non-compact convex set is proposed and  solved under the asymptotic conditions. In the end, we give a solution to the Minkowski problem for $\sigma$-finite measure on the conic domain $\Omega_C$.

 \end{abstract}

\section{Introduction}
  
 The classical Minkowski problem aims to find the necessary and/or sufficient conditions on a given finite Borel measure $\mu$ defined on the unit sphere $\sphere\subseteq \Rn$ such that $\mu$ is the surface area measure of a convex body (i.e., a convex and compact subset of $\Rn$ with nonempty interior). Its $L_p$ extension, namely the $L_p$ Minkowski problem \cite{Lut1993},  has been a central object of interest in convex geometric analysis for decades and has received extensive considerations (see e.g., \cite{CW, CWXJ, HuMaShen,  HugLYZ, LYZ2004, Umanskiy,  ZhuGX2015a,  ZhuGX2015b, ZhuGX2017}). Both the classical and $L_p$ Minkowski problems are related to  function $\varphi=t^p$ for $0\neq p\in \R$. There are versions of Minkowski problems related to other functions, for instance, the $L_0$ Minkowski or logarithmic Minkowski problems \cite{BorHegZhu, bor2013-2, Stancu2002, Stancu2003, Stancu2008,  ZhuGX2014} and the Orlicz-Minkowski problem \cite{HLYZ, huang}. 

The same problem can be asked for complete, non-compact, convex sets, where the Gauss map of their boundary is an open convex subset contained in some hemisphere (without loss of generality we set the hemisphere to be $S^n_-:=\{x\in S^n:x_{n+1}<0\}$). Then the problem can be given by:
\begin{prob}\label{p1.1}
Under what conditions on a given finite Borel measure $\mu$ defined on a convex domain $D$ in $S^n_-$, one can find a complete, non-compact, convex sets whose surface area measure is equal to $\mu$?
\end{prob} 
When restricting to convex polygon, this has been solved in the seminal papers by Alexandrov \cite{Aleks1938}. A partial solution of this problem by the Mong\`e-Amp\'ere equation on half sphere was given by Chow and Wang in their groundbreaking paper \cite{CWXJ95}. The $L_p$ Minkowski problem for all $p>0$ was recently solved by Huang and Liu in their remarkable paper \cite{Huang2021}, where the existence for the finite Borel measure $\mu$ being the $L_p$ surface area measure of a complete non-compact convex set were provided.

In view of the classical Minkowski problem and its various extensions, it is important to investigate the first variation of the volume of the convex bodies. More precisely, Khovanski\u{i}-Timorinwe \cite{KT2014}, Milman-Rotem \cite{Milman2014}, Schneider \cite{Schneider18}, and Yang-Ye-Zhu \cite{YYZ} built up the Brunn-Minkowksi theory for C-coconvex sets to solve Problem \ref{p1.1}. The set $A\subseteq C$ is called a C-coconvex set if $A$ has a finite Lebesgue measure and $A^{\bullet}=C\backslash A$ is closed and convex, where $C$ is a pointed closed convex cone with vertex at the origin and non-empty interior. To develop the Minkowski theory, Schneider \cite{Schneider18} defined co-sum of two C-coconvex sets $A_1$ and $A_2$ by
$$
A_1\oplus A_2=C\backslash(A^{\bullet}_1+A^{\bullet}_2),
$$
where $+$ is the usual Minkowski sum. Thus, Schneider \cite{Schneider18} create the following complemented Brunn-Minkowski inequality,
\begin{equation}\label{e1.1}
V_n((1-\lambda)A_1\oplus\lambda A_2)^{\frac{1}{n}}\leq (1-\lambda)V_n(A_1)^{\frac{1}{n}}+\lambda V_n(A_1)^\frac{1}{n}
\end{equation}

When the cone is studied, it is natural to admit other convex set $A$ containing some convex set $K$ such that $A\backslash K$ has finite volume. And the inequality \ref{e1.1} and equality condition could be obtained similarly from the result in Schneider's \cite{Schneider18}. But in order to specify the relation between $A$ and $K$, we would introduce the asymptotic boundary set of the non-compact convex set.

We set $\mathcal{C}_0$ to be the class of all the non-compact convex set contained in the $n+1$-dimensional upper half space $\mathbb{H}_{n+1}^+$ with the origin on the boundary and non-empty interior. For any $K\in\mathcal{C}_0$ and $t\geq 0$, we define the {\em $t$-asymptotic boundary set} of $K$ to be the set
$$
\mathsf{A}_t(K):=\left\{0\leq \langle x,e_{n+1} \rangle\leq t\right\}\bigcap \left(\bigcap_{x\in \partial K \cap (e^\perp_{n+1}+te_{n+1})}H_x^-\right),
$$
where $H_x$ is the supporting hyperplane of $K$ at $x$ and $H_x^-$ is the half space containing $K$. Then the {\em asymptotic boundary set} of $K$ to be the set
$$
\mathsf{A}(K):=\overline{\bigcup_{t\geq 0}\mathsf{A}_t(K)}. 
$$

It is clear that $\mathsf{A}(K)$ might be the whole upper half space; therefore, we define a subclass of $\mathcal{C}_0$ to be the one with compact $S(K)\cap e^{\perp}_{n+1}$, denoted by $\mathcal{C}_{c}$. In section \ref{section:SC}, we will introduce some properties of this special class. Under this special class, we could study the Brunn-Minkowski inequality for the coconvex sets $K^{c}:=\mathsf{A}(K)\backslash K$ with finite volume. In this case, we denote
$$
\mathcal{C}_b(D):=\{K\in \mathcal{C}_c(D): V_n(K^c)<\infty\}.
$$
\par
\bt \label{B-M thm}
Let $K_0,K_1\in \mathcal{C}_{b}$ with $\mathsf{A}(K_0)=a\mathsf{A}(K_1)$ for some $a>0$, and let $\lambda\in (0,1)$. Then
\begin{equation}\label{BM ine}
	V_n((1-\lambda)K_0^c\oplus\lambda K_1^c)^{\frac{1}{n}}\leq (1-\lambda)V_n(K_0^c)^{\frac{1}{n}}+\lambda V_n(K_1^c)^\frac{1}{n},
\end{equation}
where $K_0^c\oplus K_1^c=\mathsf{A}(K_0+K_1)\backslash (K_0+K_1)$. Equality holds if and only if $K_0=aK_1$.
\et

The first application of the equality condition in Theorem \ref{B-M thm} is the Minkowski uniqueness theorem for coconvex sets.
\bt\label{BM thm unique}
Let $K,L\in \mathcal{C}_b(D)$ with $\mathsf{A}(K)=\mathsf{A}(L)$. Then $K=L$, if the following identity holds on any compact set $\omega\subseteq D$:
\be\label{equ4.16}
dS(K,\cdot)=dS(L,\cdot).
\ee
\et
The fact the coconvex set $K^c$ has finite volume is crucial in the proof. Beside the uniqueness of Minkowski theorem, the existence is another essential part. First, we can get a sufficient condition since the limit of the surface area measure will go to infinity when $u\to \partial D$.
\bt\label{M thm}
Let $D$ be a convex domain in $S^n_-$ containing $-e_{n+1}$ and $f:\partial D\to \mathbb{R}_+$ being a positive continuous function on $\partial D$. Given a compact set $\omega\subseteq D$ containing $-e_{n+1}$ in its interior and a nonzero finite Borel measure $\mu$ on $D$ whose support is concentrated on $\omega$, there exists a non-compact convex set $K\in\mathcal{C}_b(D)$ with 
$$
\mathsf{A}(K)=\mathbb{H}_{n+1}^+ \cap\left[\mathop{\bigcap}\limits_{u\in \partial D} H^-(u,f(u))\right]
$$
such that 
\be
d\mu = cdS(K,\cdot) \mbox{ in } \omega,
\ee
where
$$
c=\frac{1}{nV_n(K^c)}\int_\omega h(K, u)\, d\mu(u).
$$
\et
From the nonzero finite measure with compact support on $D$, one can ask the similar questions for the nonzero $\sigma$-finite measure on $D$ since the  $\sigma$-finite measure can generate the nonzero finite measure on some compact set on $D$. However, the integral representation for the volume of coconvex set is unknown. Here, we can only answer the above question in C-coconvex set.
\bt\label{GM thm}
Let $C$ be a pointed closed convex cone with apex $o$ and $\Omega_C:=S^n\cap C^\circ$. Given a nonzero $\sigma$-finite Borel measure $\mu$ on $\Omega_C$, there exists an unique C-close set $M\in\mathcal{C}_b(\Omega_C)$ with $V_n(C\backslash M)=1$ such that for any compact $\omega\subseteq D$,
\be
d\mu = cd\bar{S}(K,\cdot) \mbox{ in } \omega,
\ee
with
$$
c=\frac{1}{n}\int_{D}\bar{h}_M\,d\mu.
$$
\et
 This answers the Minkowski problem for C-closed sets (see \cite{Schneider18}), but it can not be generated to our case.
\section{Background and Notations}\label{section 2}

Throughout this paper, $n\geq 2$ is a natural number. The notations and definations in this paper mainly follows those in \cite{SchneiderBook} for consistence.

The usual Euclidean norm is written by $\|\cdot\|$ and the origin of $\Rn$ is denoted by $o$.  Let $\{e_1,\cdots, e_n\}$ be the standard orthonormal basis of $\Rn$. 
Define $a K=\{a x: x\in K\}$ for $a\in \R$ and $K\subseteq \Rn$. For a finite measurable set $K\in \Rn$,  $V_n(K)$ refers to the Lebesgue measure of $K$ and $\cH^{n-1}$ refers to the $n-1$ dimensional Hausdorff measure. 

A subset $K\subseteq\Rn$ is convex if $(1-\lambda) x+\lambda y\in K$ for any $x, y\in K$ and $\lambda\in [0, 1]$. A subset $D\subseteq S^n_-$ is convex if any $S^n_-$-geodesic curve connecting any two points $\eta,\zeta\in D$ is in $D$ itself, equivalently the cone $\{x\in \R^{n+1}: \frac{x}{\|x\|}\in D\}\cup o$ is convex. For a set $E\subseteq \R^n$, define $conv(E)$, the convex hull of $E$, to be the smallest convex set containing $E$.  

For $\theta\in \sphere$ and $t\in \R$, let $H(\theta,t)=\{x\in \Rn: \langle x, \theta\rangle=t\}$,  $H^+(\theta,t)=\{x\in \Rn: \langle x, \theta\rangle\geq t\}$, and $H^-(\theta,t)=\{x\in \Rn: \langle x, \theta\rangle\leq t\}$.

The support function of a convex compact set $K$ containing the origin is the function $h_K: S^{n-1}\rightarrow [0, \infty]$ defined by 
$$
   h_K(u)=\max_{y\in K}\langle y, u\rangle,
$$ 
where $\langle \cdot, \cdot\rangle$ is the usual inner product on $\Rn$.  Hereafter, $S^{n-1}$ is the unit sphere of $\Rn$ which consists of all unit vectors in $\Rn$. Note that the support function $h_K$  can be extended to $\Rn\setminus \{o\}$ by 
$$
   h_K(x)=h_K(ru)=rh_K(u)
$$
for $x=ru$ with $u\in S^{n-1}$ and $r\geq 0$. Clearly $h_K: \sphere \rightarrow \R$ is sublinear.

Recall that the surface area measure $S(K, \cdot)$ of $K\in \mathcal{C}_0$ has the following geometric interpretation  (see e.g., \cite[page 111]{SchneiderBook}): for any Borel set $\Sigma \subseteq D\subseteq S_-^{n-1}$, 
\begin{equation}\label{Def:surface area----123} 
   S(K, \Sigma )=\cH^{n-1}\{x\in \partial K:  \mathrm{g}(x)\in \Sigma\},
\end{equation}  
 where $\mathrm{g}: \partial K \rightarrow D$ is the (single-valued) Gauss map of $K$, that is, $\mathrm{g}(x)\in S^{n-1}$ is the unit outer normal vector of $\partial K$ at almost everywhere $x\in \partial K$ with respect to the $(n-1)$-dimensional Hausdorff measure of $\partial K$.

 \subsection{The asymptotic boundary set of non-compact convex sets} \label{section:SC}
We start with the following definitions. Let $\mathcal{C}_0$ be the class of all the non-compact convex set contained in the $n+1$-dimensional upper half space $\mathbb{H}_{n+1}^+$ with the origin on the boundary and non-empty interior. Recall $ H(e_{n+1},t)=\{x\in \R^{n+1}: \langle x, e_{n+1}\rangle=t\}$.
\bd
For any $K\in\mathcal{C}_0$ and $t\geq 0$, we define the $t$-asymptotic boundary set of $K$ to be the set
$$
\mathsf{A}_t(K):=\left\{0\leq \langle x,e_{n+1} \rangle\leq t\right\}\bigcap \left(\bigcap_{x\in \partial K \cap  H(e_{n+1},t)}H_x^-\right),
$$
where $H_x$ is the supporting hyperplane of $K$ at $x$ and $H_x^-$ is the half space containing $K$. 
\ed
From the definition, it is clear that $K\cap \left\{0\leq \langle x,e_{n+1} \rangle\leq t\right\}\subseteq \mathsf{A}_t(K)$ and $\mathsf{A}_{t_1}\subseteq \mathsf{A}_{t_2}$ for any $0< t_1\leq t_2$. Therefore, we can consider the closure of the union of $\mathsf{A}_t(K)$ and define the following sets.
\bd\label{Dsup}
For any $K\in\mathcal{C}_0$, we define the asymptotic boundary set of $K$ to be the set
$$
\mathsf{A}(K):=\overline{\bigcup_{t> 0}\mathsf{A}_t(K)}. 
$$
\ed
Note that in general $\mathsf{A}(K)$ may be the whole upper half space, e.g. 
$$
K=\left\{(x_1,\cdots,x_{n+1})\in \Rn: x_{n+1}\geq \sum_{i=1}^n x^2_i\right\}.
$$
Moreover, the Minkowski problem may not have a unique solution in this class, see Chow-Wang \cite{CWXJ95} for more detail. Therefore, we consider the following subclass
$$
\mathcal{C}_c:=\{K\in \mathcal{C}_0: \mathsf{A}(K)\cap H(e_{n+1},0)\mbox{ is compact}\}.
$$
This subclass is exactly the collection of the non-compact convex sets with an asymptotic boundary. And moreover, we set $\mathcal{C}_c(D)$ to be the class of non-compact convex sets with support functions defined on $D$. It is clear that $D$ is convex set in $S^n_-$. Indeed the set of points where the support function of a nonempty convex set is finite is a convex cone. In this special class, we can define the coconvex sets.
\bd
Let $K\in \mathcal{C}_c$. We say that $K^c:=\mathsf{A}(K)\backslash K$ is coconvex set, if $K^c$ has positive finite Lebesgue measure. 
\ed
{\flushleft To ensure that $K^c$ is coconvex set, we define the following set} 
$$
\mathcal{C}_b(D):=\{K\in \mathcal{C}_c(D): V_n(K^c)<\infty\}.
$$
\par
Now we need to show the asymptotic boundary set are closed under Minkowski addition.
\bl \label{lem2.1}
Let  $K_0,K_1\in \mathcal{C}_{c}$ and $\lambda\in (0,1)$. Then
$$
\mathsf{A}((1-\lambda)K_0+\lambda K_1)=(1-\lambda)\mathsf{A}(K_0)+\lambda \mathsf{A}(K_1).
$$
\el
\begin{proof}
	On one hand, for any $z\in \partial((1-\lambda)K_0+\lambda K_1)\cap H(e_{n+1},t)$, there exist $x\in \partial K_0$ and $y\in \partial K_1$ such that
	$$
	z=(1-\lambda) x+\lambda y,
	$$
    which gives
    $$
    H_z^-\subseteq (1-\lambda) H_x^-+\lambda H_y^-.
    $$
    Together with the fact $\partial((1-\lambda)K_0+\lambda K_1)\cap  H(e_{n+1},t)$ is compact, we have
    \begin{align*}
    &\mathsf{A}_t((1-\lambda)K_0+\lambda K_1)\\
    =&\left\{0\leq \langle x,e_{n+1} \rangle\leq t\right\}\bigcap \left(\bigcap_{z\in \partial ((1-\lambda)K_0+\lambda K_1) \cap H(e_{n+1},t)}H_z^-\right)\\
    \subseteq & \left\{0\leq \langle x,e_{n+1} \rangle\leq t\right\}\bigcap \left(\bigcap_{z=(1-\lambda) x+\lambda y}(1-\lambda) H_x^-+\lambda H_y^-\right)\\
    \subseteq &(1-\lambda)\mathsf{A}_{t_1}(K_0)+\lambda \mathsf{A}_{t_2}(K_1)\\
    \subseteq & (1-\lambda)\mathsf{A}(K_0)+\lambda \mathsf{A}(K_1),
    \end{align*}
    where $(t_1,t_2)=\{(\max\langle x,e_{n+1}\rangle,\max\langle y,e_{n+1} \rangle): x\in \partial K_0,\ y\in \partial K_1,\ z\in\partial ((1-\lambda)K_0+\lambda K_1) \cap H(e_{n+1},t) \mbox{ and }z=(1-\lambda) x+\lambda y\}$.
    Therefore, by definition \ref{Dsup}, we have 
    $$
    \mathsf{A}((1-\lambda)K_0+\lambda K_1)\subseteq (1-\lambda)\mathsf{A}(K_0)+\lambda \mathsf{A}(K_1).
    $$
    
    On the other hand, it is clear that
    \begin{align*}
    	& (1-\lambda)\mathsf{A}_{t_1}(K_0)+\lambda \mathsf{A}_{t_2}(K_1)\\
    	\subseteq & \mathsf{A}_{t_0}((1-\lambda)K_0+\lambda K_1)\\
    	\subseteq & \mathsf{A}((1-\lambda)K_0+\lambda K_1),
    \end{align*}
    where $t_0=\max\langle (1-\lambda) x+\lambda y, e_{n+1}\rangle: x\in \partial K_0\cap H(e_{n+1},t)\mbox{ and }y\in \partial K_1 \cap H(e_{n+1},t)$. Therefore, by definition \ref{Dsup}, we have 
    $$
    \mathsf{A}((1-\lambda)K_0+\lambda K_1)\supseteq (1-\lambda)\mathsf{A}(K_0)+\lambda \mathsf{A}(K_1).
    $$
\end{proof}

\subsection{Wulff shapes}\label{section:WS}
Let $C^+(D)$ be the class of non-negative functions on $D\subseteq S^n_-$ containing $-e_{n+1}$. Associated to each $f\in  C^{+}(\bar{D})$ with $f(-e_{n+1})=0$, one can define a convex set $K_{f,D}\in \mathcal{C}_c$ (see Wulff shapes in \cite[Section 7.5]{SchneiderBook}) by  
$$
K_{f,D}=\mathop{\bigcap}\limits_{u\in D} H^-(u,f(u)).
$$ 
The convex set $K_{f,D}$ is called the {\em Wulff shape}, or Aleksandrov body associated to $f \in C^{+}(\bar{D}).$ The Aleksandrov body provides a powerful tool in convex geometry and plays crucial roles in this paper. Here we list some important properties for  the Aleksandrov body which will be used in later context. These properties and the proofs can be found in section 7.5 in \cite{SchneiderBook}.  First of all,  if $f\in C^+(\bar{D})$ is the support function of a convex body $K\in \mathcal{C}_c$, then $K=K_{f,D}$. Secondly, for $f\in C^+(\bar{D})$,  $h_{K_{f,D}}(u)\leq f(u)$ for all $u\in D$, but the fact $h_{K_{f,D}}(u)=f(u)$ almost everywhere with respect to $S(K_{f,D}, \omega)$ may be false if $S(K_{f,D}, \omega)$ is unbounded. In this section, $\omega$ is always some nonempty compact set $\omega\subseteq D$ with a positive distance away from the boundary of $D$ (see Section 8 in \cite{Schneider18}). Therefore, we need some proper restrictions on $D$.

\bd
Give a convex domain $D\subseteq S^n_-$ containing $-e_{n+1}$ and a initial fixed continuous function $f$ on $\bar{D}$ with
$$
f(u)
\begin{cases}
>0 & u\in \bar{D}\backslash\{-e_{n+1}\}\\
=0& u=-e_{n+1},	
\end{cases}
$$
we call $(D,f)$ is irreducible if
$$
\mathsf{A}(K_{f,D})=\mathbb{H}_{n+1}^+ \cap\left[\mathop{\bigcap}\limits_{u\in \partial D} H^-(u,f(u))\right]=:K_{f,\partial D}.
$$
\ed

Since one concerns a slight difference from coconvex Wulff shapes, we reprove some results in \cite[Section 8]{Schneider18}.

\bl\label{lem2.2}
Give a convex domain $D\subseteq S^n_-$ containing $-e_{n+1}$ and a initial fixed continuous function $f$ on $\bar{D}$ with 
$$
f(u)
\begin{cases}
	>0 & u\in \bar{D}\backslash\{-e_{n+1}\}\\
	=0& u=-e_{n+1},	
\end{cases}
$$
if $(D,f)$ is irreducible, there exists a unique noncompact convex surface $K$ in $\mathcal{C}_c$ with support functions $h_K$ on $\bar{D}$ such that $h_K=f$ almost everywhere with respect to $S(K, \omega)$ for any compact set $\omega\subseteq D$.
\el
\begin{proof}
We set 
$$
K_{f,D}=\bigcap_{u\in D}H^-(u,f(u)).
$$
By the Wulff shapes introduced by \cite[Section 7.5]{SchneiderBook}, for a $t>0$, there exists a unique convex body 
$$
K_t:=K_{f,D}\bigcap H^-(e_{n+1},t).
$$
It is clear that $K_{t_1}\subseteq K_{t_2}$ for $0<t_1<t_2$, and 
$$
\lim_{t\to\infty}K_t=K_{f,D}.
$$
Moreover, for any compact $\omega\subseteq D$ away from the boundary of $D$, there exists a $\varepsilon>0$ such that
$$
\|\theta-u\|\geq \varepsilon\ \mbox{for }\theta\in \partial D\mbox{ and }u\in \omega
$$
Since $(D,f)$ is irreducible, for any $x\in \mathsf{A}(K_{f,D})$ and $u\in \omega$, there exists a $\theta\in \partial D$ and $a_0>0$ such that
$$
\langle x-\langle x,\theta\rangle\theta,\theta \rangle=0
$$
and
\begin{equation}\label{equ2.4}
	\langle x-\langle x,\theta\rangle\theta,u \rangle\leq -a_0\|x\|.
\end{equation}
Indeed, if $\theta\notin \partial S_-^n$ see \cite[Section 8]{Schneider18}; if $\theta\in \partial S_-^n$, the results above follows the fact $\frac{x}{\|x\|}= e_{n+1}$ and $\omega$ being compact. 
\par
Now we define a convex set
$$
K_{\omega}=\mathsf{A}(K_{f,D})\cap \left[\bigcap_{u\in \omega}H^-(u,f(u))\right].
$$
Note that if for some $u\in \omega$ with $x\notin H^-(u,h_{K_{\omega}}(u))$ and Equation \ref{equ2.4}, we can obtain
$$
\|x\|\leq \frac{1}{a_0}\left(\langle x,\theta\rangle\langle\theta,u \rangle-\langle x,u \rangle\right)\leq \frac{1}{a_0}\left(\langle\theta,u \rangle h_{K_\omega}(\theta)-h_{K_\omega}(u)\right).
$$
Hence $\mathsf{A}(K_{f,D})\backslash K_{\omega}$ is bounded since $\partial D$ and $\omega$ are compact and $h_{K_\omega}$ is continuous.
\par
As a result, there exists $t>0$ such that $S(K_t,\omega)=S(K,\omega)$.
Therefore, by \cite[Section 7.5]{SchneiderBook}, $h_K=f$ almost everywhere with respect to $S(K_t, \omega)=S(K,\omega)$.
\end{proof}
In special case, when only defining support function on the boundary of some convex set $D\subseteq S_-^n$, we can deduce a special convex surface in $\mathcal{C}_c$.

\br
This special class may not be homothetic after some large $t>0$. We could consider a lune area $D$ on $S_-^n$ with $\theta_0$ and $-\theta_0$ being the vertices, and set
$$
K=\mathbb{H}_{n+1}^+\cap H^-(\theta_0,1)\cap H^-(-\theta_0,1)\cap \left[\bigcap_{\theta\in \partial D\backslash \{\theta_0,-\theta_0\}}H(\theta,-\langle e_{n+1},\theta\rangle)\right].
$$
The convex set $K$ is generated by a convex domain in $S_-^n$, but for any two different sufficiently large $t_1,t_2>0$, the parallel sections $K\cap H^+(e_{n+1},t_1)$ and $K\cap H^+(e_{n+1},t_2)$ are not homothetic.
\er

\section{The Brunn-Minkowski inequality} \label{Section:OM-Inequality}

This section is dedicated to establish the Brunn-Minkowski inequality. The following proof of Theorem \ref{B-M thm} comes from the methods in \cite[Section 3]{Schneider18}. 

\begin{proof}[proof of Theorem \ref{B-M thm}]
	Let $K$ and $L$ be the convex sets in Theorem \ref{B-M thm}. Without loss of generality, for $\nu\in \{0,1\}$ we set 
	$$
	V_n(K^c_\nu)=1.
	$$
	Otherwise, setting
	$$
	\bar{K}_\nu:=[V_n(K^c_\nu)]^{-\frac{1}{n}}K_\nu
	$$
	and
	$$
	\bar{\lambda}:=\frac{\lambda [V_n(K^c_1)]^{\frac{1}{n}}}{(1-\lambda) [V_n(K^c_0)]^{\frac{1}{n}}+\lambda [V_n(K^c_1)]^{\frac{1}{n}}}.
	$$
	Then $V_n(S(\bar{K}^c_\nu)=1$ and $V_n((1-\bar{\lambda})\bar{K}^c_0\oplus \bar{\lambda}\bar{K}^c_1)\leq 1$ give the result.
	
	Now for $t>0$, we recall 
	$$
	H^-(e_{n+1},t):=\{x\in \R^{n+1}:\langle x,e_{n+1}\rangle\leq t\}
	$$
	and set
	$$
	v_\nu(t):=V_{n-1}(K_\nu\cap H(e_{n+1},t))\mbox{ and }w_\nu(t):=V_n(K_\nu\cap H^-(e_{n+1},t))
	$$
	for $\nu\in \{0,1\}$. Then
	$$
	w_\nu(t)=\int_0^t v_\nu(s)\,ds
	$$
	and
	$$
	w'_\nu(t)=v_\nu(t)>0\mbox{ for }0<t<\infty.
	$$
	Now set $z_\nu$ being the inverse function of $w_\nu$, then
	$$
	z'_\nu(\tau)=\frac{1}{v_\nu(z_\nu(\tau))}\mbox{ for }0<\tau<\infty.
	$$
	Define
	$$
	\begin{cases}
		D_\nu(\tau):=K_\nu\cap H(e_{n+1},z_\nu(\tau)),\\
		K_\lambda:=(1-\lambda)K_0+\lambda K_1,\\
		z_\lambda(\tau):=(1-\lambda)z_0(\tau)+\lambda z_1(\tau).
	\end{cases}
	$$
	Then it is clear that
	$$
	K_\lambda\cap H(e_{n+1},z_\lambda(\tau)) \supseteq (1-\lambda)D_0(\tau)+\lambda D_1(\tau).
	$$
	For $\tau>0$ we define
	$$
	b_\nu(\tau):=V_n(K^c_\nu\cap H^-(e_{n+1},z_\nu(\tau)))=V_n(\mathsf{A}(K_\nu)\cap H^-(e_{n+1},z_\nu(\tau)))-\tau
	$$
	and $f(\tau):=V_n(K_\lambda\cap (e^\perp_{n+1})^-_{z_\nu(\tau)})$, then
	\be\label{e3.4}
	V_n(K^c_\lambda \cap H^-(e_{n+1},z_\nu(\tau)))=V_n(\mathsf{A}(K_\lambda)\cap H^-(e_{n+1},z_\nu(\tau)))-f(\tau),
	\ee
	where
	\begin{align*}
		f(\tau)&=\int_0^{z_\lambda(\tau)}V_{n-1}(K_\lambda\cap H(e_{n+1},\zeta))\,d\zeta\\
		&=\int_0^\tau V_{n-1}(K_\lambda\cap H(e_{n+1},z_\nu(t)))z'_\lambda(t)\,dt\\
		&\geq \int_0^\tau V_{n-1}((1-\lambda)D_0(t)+\lambda D_1(t))z'_\lambda(t)\,dt\\
		&\geq \int_0^\tau \left[(1-\lambda)v_0(z_0(t))^{\frac{1}{n-1}}+\lambda v_1(z_1(t))^{\frac{1}{n-1}}\right]^{n-1}\left[\frac{1-\lambda}{v_0(z_0(t))}+\frac{\lambda}{v_1(z_1(t))}\right]\,dt\\
		&\geq \tau
	\end{align*}
\par
    Now for $t>0$ we define
    $$
    g_\nu(t):=V_n(\mathsf{A}(K_\nu)\cap H^-(e_{n+1},t))\mbox{ and }g_\lambda(t):=V_n(\mathsf{A}(K_\lambda)\cap H^-(e_{n+1},t))
    $$
    and by Lemma \ref{lem2.1}, without loss of generality, we set $a\leq 1$, then it is clear that $g_1(t)=a^{-n}g_0(at)$ and $g_\lambda(t)=\beta^{-n}g_0(\beta t)$, where $\beta= \frac{a}{(1-\lambda)a+\lambda}$. Note that by Brunn theorem, $g_0^{\frac{1}{n}}(t)$ is concave and increasing function. Hence using the mean value theorem, there exists $\eta_0,\eta_1\geq 0$ such that
    $$
    g_0^{\frac{1}{n}}(\beta z_\lambda(\tau))-g_0^{\frac{1}{n}}(z_0(\tau))=\eta_0(\beta z_\lambda(\tau)-z_0(\tau))=\eta_0\tilde{\beta}(az_1(\tau)-z_0(\tau))
    $$
    and
    $$
    g_0^{\frac{1}{n}}(a z_1(\tau))-g_0^{\frac{1}{n}}(\beta z_\lambda(\tau))=\eta_1(a z_1(\tau)-\beta z_\lambda(\tau))=\eta_1(1-\tilde{\beta})(az_1(\tau)-z_0(\tau))
    $$
    where $\tilde{\beta}=\frac{\lambda}{(1-\lambda)a+\lambda}$. Therefore, we obtain
    $$
    \eta_1(1-\tilde{\beta})\left(g_0^{\frac{1}{n}}(\beta z_\lambda(\tau))-g_0^{\frac{1}{n}}(z_0(\tau))\right)=\eta_0\tilde{\beta}\left(g_0^{\frac{1}{n}}(a z_1(\tau))-g_0^{\frac{1}{n}}(\beta z_\lambda(\tau)) \right)
    $$
    which implies
    $$
    \eta_1(1-\tilde{\beta})\left(\beta g_{\lambda}^{\frac{1}{n}}(z_\lambda(\tau))-g_0^{\frac{1}{n}}(z_0(\tau))\right)=\eta_0\tilde{\beta}\left(ag_1^{\frac{1}{n}}(z_1(\tau))-\beta g_{\lambda}^{\frac{1}{n}}(z_\lambda(\tau)). \right)
    $$
    Hence
    \be \label{e3.6}
    g_{\lambda}^{\frac{1}{n}}(z_\lambda(\tau))=(1-\tilde{\lambda})g_0^{\frac{1}{n}}(z_0(\tau))+\tilde{\lambda} a g_1^{\frac{1}{n}}(z_1(\tau))
    \ee
    where $\tilde{\lambda}=\frac{\eta_0\lambda}{\eta_0\lambda+\eta_1(1-\lambda)a}$. Here, $\tilde{\lambda}\in [0,1]$ for any fixed $\lambda\in [0,1]$ and $a>1$. By Equation \ref{e3.4} and \ref{e3.6}, we have
	\begin{align*}
		& V_n(K^c_\lambda \cap H^-(e_{n+1}, z_\lambda(\tau)))\\
		=&V_n(\mathsf{A}(K_\lambda)\cap H^-(e_{n+1}, z_\lambda(\tau)))-f(\tau)\\
		=&g_{\lambda}(z_\lambda(\tau))-f(\tau)\\
		=&\left[(1-\tilde{\lambda})g_0^{\frac{1}{n}}(z_0(\tau))+\tilde{\lambda} a g_1^{\frac{1}{n}}(z_1(\tau))\right]^n-f(\tau)\\
		\leq & \left[(1-\tilde{\lambda})g_0^{\frac{1}{n}}(z_0(\tau))+\tilde{\lambda} g_1^{\frac{1}{n}}(z_1(\tau))\right]^n-f(\tau)\\
		=&\left[(1-\tilde{\lambda})(b_0(\tau)+\tau)^{\frac{1}{k}}+\tilde{\lambda} (b_1(\tau)+\tau)^{\frac{1}{k}}\right]^k-f(\tau)\\
		=&\left[(b_0(\tau)+\tau)^{\frac{1}{k}}+\tilde{\lambda} (b_1(\tau)-b_0(\tau))h(\tau)\right]^k-f(\tau)\\
		=&b_0(\tau)+\tau -f(\tau)+\sum_{r=1}^k\left(\begin{matrix}
			k\\r
		\end{matrix}\right)(b_0(\tau)+\tau)^{\frac{k-r}{k}}\left[\tilde{\lambda} (b_1(\tau)-b_0(\tau))h(\tau)\right]^r.
	\end{align*}
	Together with $b_0(\tau)\to 1$, $f(\tau)\geq \tau$, $h(\tau):=\frac{1}{k}(b(\tau)+\tau)^{\frac{1}{k}-1}$, $b(\tau)$ is a function between $b_0(\tau)$ and $b_1(\tau)$, and $b_1(\tau)-b_0(\tau)\to 0$, we can conclude that
	\be\label{e3.7}
	V_n(K^c_\lambda)\leq 1.
	\ee
	\par
	Suppose the Equation \ref{e3.7} holds for $K_0$ and $K_1$. Then, $f(\tau)=\tau$ for all $\tau>0$ and $a=1$, which give
	$$
	K_\lambda\cap H^-(e_{n+1}, z_\lambda(\tau))=(1-\lambda)	K_0\cap H^-(e_{n+1}, z_0(\tau))+\lambda K_1\cap H^-(e_{n+1}, z_1(\tau)).
	$$
	Thus, by the same argument in \cite[Section 3]{Schneider18}, $K_0=K_1$.
\end{proof}

\section{The related mixed volume}
In this section, we develop the mixed volume of coconvex sets under the same asymptotic boundary set. However, there does not exist an integral representation of the volume of the coconvex sets.
\par
Let $K\in \mathcal{C}_c(D)$. Since $o\in \partial K$, The support function $h_K$ is a finite non-negative function defined on $\bar{D}$ with $D\subseteq S_-^n$. The surface area measure $S(K,\omega)$ of $K$ is defined by
$$
S(K,\omega):=\mathcal{H}^{n-1}(\mathsf{g}^{-1}(K,\omega))
$$ 
for Borel set $\omega\in D$. Recall $\mathsf{g}^{-1}(K,\omega)$ is the inverse Gauss image of $K$ where the outer normal vector falling in $\omega$. Since the integral over the surface of the non-compact convex set is unbounded, we restrict the integral on a compact closed set $\omega\subseteq D$ with a positive distance away from the boundary of $D$ ($\omega$ always represents this set here and below in this section). Recall the wulff shape associated with $(K,\omega)$ to be a closed convex set
$$
K_\omega:=A(K)\cap  \left[\bigcap_{u\in \omega} H^-(u, h_K(u))\right].
$$
By Lemma \ref{lem2.2}, $A(K)\backslash K_\omega$ is bounded, hence the surface integral over $\omega$ is bounded. First we need to show the continuity of surface area measure on $\omega$.
\bd
Let $\{K_i\}_{i\geq 1}\subseteq \mathcal{C}_c$ and $K\in \mathcal{C}_c$ with $\mathsf{A}(K_i)=\mathsf{A}(K)$ for any $i$. We write $K_i\to K$ if
$$
d_H(\mathsf{A}(K)\cap K_i,\mathsf{A}(K)\cap K)\to 0,
$$
where $d_H(\cdot,\cdot)$ is the Hausdorff distance.
\ed
Here, since $\{K_i\}_{i\geq 1}\subseteq \mathcal{C}_c$ and $K\in \mathcal{C}_c$ with $\mathsf{A}(K_i)=\mathsf{A}(K)$ for any $i$, the domains of the support function of $\{K_i\}_{i\geq 1}$ are same as $K$, denoted by $D$ again. Therefore, for any closed compact  set $\omega\subseteq D$, we have the following weak convergence of surface area measures defined on $\omega$.
\bl
Let $\{K_i\}_{i\geq 1}\subseteq \mathcal{C}_c$ and $K\in \mathcal{C}_c$ with $\mathsf{A}(K_i)=\mathsf{A}(K)$ for any $i$. If $K_i\to K$, then $S(K_i,\omega)\to S(K,\omega)$ weakly on $\omega$. That is, for any continuous function $f:\omega\to \mathbb{R}$, we have
$$
\int_\omega f(u)\,dS(K_i,u)\to \int_\omega f(u)\,dS(K,u)
$$
\el
\bpf
Suppose that $K_i\to K$ and set 
$$
(\mathsf{A}(K))_t=\mathsf{A}(K)\cap H^-(e_{n+1},t)
$$
For sufficiently large $t>0$, we have
$$
d_H\left((K_i)_\omega\cap (\mathsf{A}(K))_t, K_\omega\cap (\mathsf{A}(K))_t\right)\to 0.
$$
Note that by Lemma \ref{lem2.2}, for all $i$, the surface area measure $S(K_i,\omega)$ and $S(K,\omega)$ are equivalent to $S((K_i)_{\tilde{\omega}},\omega)$ and $S(K_{\tilde{\omega}},\omega)$ for some compact sets $\nu$ and $\title{\omega}$ satisfying
$$
\omega\subseteq \mbox{int}(\nu)\subseteq\nu\subseteq \mbox{int}(\tilde{\omega})\subseteq \tilde{\omega}\subseteq D;
$$ 
moreover, we have
$$
S((K_i)_{\tilde{\omega}},D\backslash \tilde{\omega})=S(K_{\tilde{\omega}},D\backslash \tilde{\omega})=0.
$$
By Tietze's extension theorem, for any continuous function $f$ on $\omega$, there exists a continuous function $F:S^{n}\to \mathbb{R}$ such that
$$
F(u)=
\begin{cases}
	f(u) & \mbox{on }\omega\\
	0 & \mbox{on }S^n\backslash \nu.
\end{cases}
$$
Therefore, the Hausdorff convergence of convex bodies implies
$$
\int_{S^n}F(u)\,dS((K_i)_\omega\cap (\mathsf{A}(K))_t,u)\to \int_{S^n}F(u)\,dS(K_\omega\cap (\mathsf{A}(K))_t,u),
$$
which yields
$$
\int_{\nu}F(u)\,dS((K_i)_{\tilde{\omega}},u)\to \int_{\nu}F(u)\,dS(K_{\tilde{\omega}},u),
$$
By the arbitrary choice of $\nu$ and compactness, we can obtain
$$
\int_{\omega}f(u)\,dS((K_i)_{\tilde{\omega}},u)\to \int_{\omega}f(u)\,dS(K_{\tilde{\omega}},u),
$$
hence
$$
\int_{\omega}f(u)\,dS(K_i,u)\to \int_{\omega}f(u)\,dS(K,u),
$$
\epf
We are now ready to state and prove the variational formula regarding the volume of the coconvex sets. Let $K\in \mathcal{C}_c(D)$ and $f:\omega\to \mathbb{R}$ be a continuous function on the compact set $\omega\subseteq D$. Define $f_\tau: \omega\to \mathbb{R}$ with $|\tau|<\tau_0$ for some sufficiently small $t_0>0$, by
$$
f_\tau(u)=h_K(u)+\tau f(u)\mbox{ for }u\in \omega.
$$
By the continuity of $h_K(u)$ and $f(u)$ on the compact set $\omega$, it is easy to verify when
$$
0<\tau_0\leq \frac{\min_{\omega}h_K(u)}{\max_{\omega}f(u)},
$$
$f_\tau$ is also a positive and continuous function on $\omega$ and $f_\tau \to h_K(u)$ uniformly on $\omega$ as $\tau \to 0$. Hence, defining 
\be\label{equ4.8}
K_\tau:=\left[\bigcap_{u\in \omega}H^-(u,f_\tau(u))\right]\cap \mathsf{A}(K),
\ee
we have $\mathsf{A}(K_\tau)=\mathsf{A}(K)$ and $K_\tau$ converges to $K_\omega$ as $\tau\to 0$.
\bt\label{thm4.1}
Let $K\in \mathcal{C}_b(D)$ and $f:\omega\to \mathbb{R}$ be a continuous function on the compact set $\omega\subseteq D$. For $K_\tau$ being defined by (\ref{equ4.8}), one has
\be
\frac{dV_n(K_\tau^c)}{\tau}\bigg|_{\tau=0}=\lim_{\tau\to 0}\frac{V_n(K^c_\tau)-V_n(K_\omega^c)}{\tau}=-\int_\omega f(u)\,dS(K_\omega,u)
\ee
\et
\bpf
Recall $M_t:=M\cap H^-(e_{n+1},t)$ for $M\in \mathcal{C}_c$. By Lemma \ref{lem2.2}, there exists a sufficiently large $t>0$ such that $K_\tau\cap H^+(e_{n+1},t)=K_\omega\cap H^+(e_{n+1},t)$. Indeed, by Lemma \ref{lem2.2}, for any $\tau\in [-\tau_0,\tau_0]$, $A(K)\backslash (K_\tau)_\omega$ is bounded, then there exists a sufficient large $t>0$ such that 
$$
A(K)\backslash (K_\tau)_\omega\subseteq H^-(e_{n+1},t)\mbox{ for all }\tau\in [-\tau_0,\tau_0].
$$ 
Hence, we have
\be\label{equ4.10}
	\lim_{\tau\to 0}\frac{V_n(K^c_\tau)-V_n(K_\omega^c)}{\tau}	=\lim_{\tau\to 0}\frac{V_n((K_\omega)_t)-V_n((K_\tau)_t)}{\tau} 
\ee
Now we set the continuous function by
$$
F(u)=
\begin{cases}
	f(u)&\mbox{on }\omega\\
	0 & \mbox{on } S^n\backslash D.
\end{cases}
$$
By the volume difference of convex sets and Lemma \ref{lem2.2}, we have
$$
\lim_{\tau\to 0}\frac{V_n((K_\omega)_t)-V_n((K_\tau)_t)}{\tau}=-\int_{S^n} F(u)\, dS((K_\omega)_t,u)=-\int_\omega f(u)\, dS(K_\omega,u).
$$
Here, by the definition of $K_\omega$, $S(K_\omega,D\backslash \omega)=0$.
\epf
Based on the variational formula, we can define the mixed volume for coconvex sets.
\bd\label{def4.2}
Let $K,L\in \mathcal{C}_b(D)$ with $\mathsf{A}(K)=\mathsf{A}(L)$. Define the mixed volume of $K_\omega^c$ and $L_\omega^c$ on $\omega$, denoted by $V_\omega(K_\omega^c,L_\omega^c)$ as
\be
V_\omega(K_\omega^c,L_\omega^c):=\frac{1}{n}\lim_{\tau\to 0}\frac{V_n((1-\tau)K_\omega^c\oplus\tau L_\omega^c)-V_n(K_\omega^c)}{\tau}\ee
\ed
{\flushleft It is easy to check}
$$
V_\omega(K_\omega^c,L_\omega^c)=\frac{1}{n}\int_\omega (h_K(u)-h_L(u))\,dS(K_\omega,u)
$$
from Theorem \ref{thm4.1}. We now prove the following Minkowski inequality.
\bt\label{thm4.2}
Let $K,L\in \mathcal{C}_b(D)$ with $\mathsf{A}(K)=\mathsf{A}(L)$. We have
\be\label{equ4.12}
\left(V_\omega(K_\omega^c,L_\omega^c)+V_n(K_\omega^c)\right)^n\leq V_n(K_\omega^c)^{n-1}V_n(L_\omega^c)
\ee
with equality if and only if $K_\omega=L_\omega$.
\et
\bpf
By the Definition \ref{def4.2} and Theorem \ref{B-M thm}, we can show
\begin{align*}
	V_\omega(K_\omega^c,L_\omega^c)=& \frac{1}{n}\lim_{\tau\to 0}\frac{V_n((1-\tau)K_\omega^c\oplus\tau L_\omega^c)-V_n(K_\omega^c)}{\tau}\\
	\leq & \frac{1}{n}\lim_{\tau\to 0}\frac{\left[(1-\tau)\left(V_n(K_\omega^c)\right)^\frac{1}{n}+\tau \left(V_n(L_\omega^c)\right)^\frac{1}{n}\right]^n-V_n(K_\omega^c)}{\tau}\\
	=& V_n(K_\omega^c)^{\frac{n-1}{n}}V_n(L_\omega^c)^{\frac{1}{n}}-V_n(K^c_\omega),
\end{align*}
where the equality holds if and only if $K_\omega=L_\omega$. 
\epf
Since for any compact set $\omega\subseteq D$, the Mixed volume and Minkowski inequality are well-defined, we can generate it to the whole coconvex sets by a sequence of convex sets $\{\omega_i\}_{i\in \mathbb{N}}$ with $\omega_i\subseteq \omega_{i+1}$ and $\cup_i \omega_i=D$. By the definition of $K_\omega$, it is clear that $K_{\omega_{i}}\supseteq K_{\omega_{i+1}}$; hence, $K=\cap_i K_{\omega_{i}}$. Moreover, for $K\in \mathcal{C}_b$, $V_n(K^c)=\lim_{i\to \infty}V_n(K^c_{\omega_{i}})$. Then we can get the following definition and theorem from above.
\bd\label{def4.3}
Let $K,L\in \mathcal{C}_b(D)$ with $\mathsf{A}(K)=\mathsf{A}(L)$. Define the mixed volume of $K^c$ and $L^c$, denoted by $V(K^c,L^c)$ as
\be
V(K^c,L^c):=\frac{1}{n}\lim_{\tau\to 0}\frac{V_n((1-\tau)K^c\oplus\tau L^c)-V_n(K^c)}{\tau}.
\ee
\ed
The following lemma shows the definition \ref{def4.3} is well-defined.
\bl\label{lem4.2}
Let $K,L\in \mathcal{C}_b(D)$ with $\mathsf{A}(K)=\mathsf{A}(L)$. Then
$$
V(K^c,L^c)=\frac{1}{n}\int_D (h_K(u)-h_L(u))\,dS(K,u).
$$
\el
\bpf
It is sufficient to show 
\be\label{equ4.14}
\lim_{i\to \infty}\frac{1}{n}\int_{\omega_{i}} (h_K(u)-h_L(u))\,dS(K_{\omega_{i}},u)=\frac{1}{n}\int_D (h_K(u)-h_L(u))\,dS(K,u).
\ee
By Theorem \ref{thm4.1}, one has for any $\omega_{i}$
\begin{align*}
	&\left|\frac{1}{n}\int_{\omega_{i}} (h_K(u)-h_L(u))\,dS(K_{\omega_{i}},u)\right|\\
	\leq &	\left|V_n(K_{\omega_i}^c)^{\frac{n-1}{n}}V_n(L_{\omega_i}^c)^{\frac{1}{n}}-V_n(K^c_{\omega_i})\right|\\
	\leq & \left|V_n(K_{\omega_i}^c)^{\frac{n-1}{n}}V_n(L_{\omega_i}^c)^{\frac{1}{n}}\right|+\left|V_n(K^c_{\omega_i})\right|\\
	\leq & \left|V_n(K^c)^{\frac{n-1}{n}}V_n(L^c)^{\frac{1}{n}}\right|+\left|V_n(K^c)\right|,
\end{align*}
On the other hand, since $\mathsf{A}(K)=\mathsf{A}(L)$, for any compact set $\omega\subseteq D$, there exists $\omega_N$ such that $\omega\subseteq \omega_{N}$, then we can show
\begin{align*}
	&\int_{\omega} (h_K(u)-h_L(u))\,dS(K_{\omega_{i}},u)\\
	= &\int_{\omega\cap \mbox{int}(D)} (h_K(u)-h_L(u))\,dS(K_{\omega_{i}},u)+\int_{\omega\cap \partial D} (h_K(u)-h_L(u))\,dS(K_{\omega_{i}},u)\\
	= &\int_{\omega\cap \mbox{int}D} (h_K(u)-h_L(u))\,dS(K_{\omega_{i}},u)\\
	=& \int_{\omega\cap\partial\omega_N\cap \mbox{int}D} (h_K(u)-h_L(u))\,dS(K_{\omega_{i}},u)+\int_{\omega\cap\mbox{int}(\omega_N\cap D)} (h_K(u)-h_L(u))\,dS(K_{\omega_{i}},u)\\
	=& \int_{\omega\cap\partial\omega_N\cap \mbox{int}D} (h_K(u)-h_L(u))\,dS(K_{\omega_{i}},u)+\int_{\omega\cap\mbox{int}(\omega_N\cap D)} (h_K(u)-h_L(u))\,dS(K,u).
\end{align*}
Now since $\cup_i \omega_i=D$, we can choose sufficient large $i$ such that $\omega\cap\partial\omega_N\cap \mbox{int}D\subseteq \mbox{int}(\omega_{i})$, which gives
$$
\int_{\omega} (h_K(u)-h_L(u))\,dS(K_{\omega_{i}},u)=\int_{\omega} (h_K(u)-h_L(u))\,dS(K,u).
$$
Therefore, for any compact set $\omega\subseteq D$, we obtain
$$
\lim_{i\to \infty}\frac{1}{n}\int_{\omega} (h_K(u)-h_L(u))\,dS(K_{\omega_{i}},u)=\frac{1}{n}\int_{\omega} (h_K(u)-h_L(u))\,dS(K,u).
$$
By compact convergence and boundness, the equation \ref{equ4.14} is true.
\epf
Then the Minkowski inequality for the coconvex sets is followed by Theorem \ref{thm4.2} and Lemma \ref{lem4.2}.
\bt\label{thm4.3}
Let $K,L\in \mathcal{C}_b(D)$ with $\mathsf{A}(K)=\mathsf{A}(L)$. We have
\be\label{equ4.15}
\left(V(K^c,L^c)+V_n(K^c)\right)^n\leq V_n(K^c)^{n-1}V_n(L^c)
\ee
with equality if and only if $K=L$.
\et

The following result is for the uniqueness theorem in $\mathcal{C}_b$ by the surface area measure. In particular, the uniqueness of the Minowski problem for non-compact convex set in $\mathcal{C}_b$ is true.

\bpf[Proof of Theorem \ref{BM thm unique}]
By Theorem \ref{thm4.3}, one has
\be\label{equ4.17}
V(K^c,L^c)+V_n(K^c)\leq V_n(K^c)^{\frac{n-1}{n}}V_n(L^c)^{\frac{1}{n}}
\ee
and
\be\label{equ4.18}
V(L^c,K^c)+V_n(L^c)\leq V_n(K^c)^{\frac{1}{n}}V_n(L^c)^{\frac{n-1}{n}}
\ee
Adding Equations (\ref{equ4.17}) and (\ref{equ4.18}), together with the fact $V(K^c,L^c)+V(L^c,K^c)$ by Lemma \ref{lem4.2} and Equation (\ref{equ4.16}), we can obtain
$$
V_n(K^c)+V_n(L^c)\leq V_n(K^c)^{\frac{n-1}{n}}V_n(L^c)^{\frac{1}{n}}+V_n(K^c)^{\frac{1}{n}}V_n(L^c)^{\frac{n-1}{n}},
$$
which implies
$$ (V_n(K^c)^{\frac{n-1}{n}}-V_n(L^c)^{\frac{n-1}{n}})(V_n(K^c)^{\frac{1}{n}}-V_n(L^c)^{\frac{1}{n}})\leq 0,
$$
where the ``$=$'' holds only, which means
$$
K=L.
$$
\epf

\section{The Minkowski problem for non-compact convex sets}

In this section, we will first provide a solution to the Minkowski problem in $\mathcal{C}_b$. But we need to show a Lemma as the one in \cite[Section 9]{Schneider18} first.

\bl\label{lem5.1}
Let $L\in \mathcal{C}_c(D)$ with $L=\mathsf{A}(L)$ and $\omega\subseteq D$ be a compact set. Then the following facts are true.
\begin{enumerate}
	\item To every bounded set $B\subset L$ there is a number $t>0$ such that if
	$$
	H(u,\tau)\cap B\neq \emptyset \mbox{ with }u\in \omega,
	$$
	then $H(u,\tau)\cap L\subseteq L_t$.
	\item For any $K\in \mathcal{C}_c(D)$ with $\mathsf{A}(K)=L$. There is a constant $t>0$ not depending on $K$ such that
	$$
	C\cap H_t\subseteq K_\omega.
	$$
\end{enumerate}
\el
\bpf
To show 1, let $B\subseteq L$ be a bounded set. We choose a $s>0$ such that $B\subseteq \mathsf{A}(L)_s$. For any $u\in \omega$, we set $H(u,\delta)$ ($\delta<0$) to be the supporting hyperplane of $\mathsf{A}(L)_s$ such that
$$
\mathsf{A}(L)_s\subseteq H^+(u,\delta).
$$
Choose $z\in \mathsf{A}(L)_s\cap H(u,\delta)$ and $x$ to be the point in $\mathsf{A}(L)\cap H(u,\delta)$ with maximal norm. It is clear that $x\in \partial A(L)$ and the outer normal vector of $\mathsf{A}(L)$ at $x$, denoted by $v$, is on $\partial D$. Since $\omega$ is compact in $D$, there exists a $a_0>0$ such that $\sin \psi \geq a_0$ where $\psi$ is the angle between $u$ and $v$. Now we set
$$
y=[o,x]\cap H(u,s),
$$
where $y\in \mathsf{A}(L)$ by the compactness of $\mathsf{A}(L)$. The triangle with vertice $x,y,z$ and angle $\alpha$ at $z$ and $\beta$ at $x$ gives 
$$
\|x-y\|=\frac{\sin \alpha}{\sin \beta}\|y-z\|\leq \frac{\mbox{diam}(\mathsf{A}(L)_s)}{a_0}.
$$
Here, by convexity of $\mathsf{A}(L)_s$, we have $\sin \beta\geq \sin \psi \geq a_0$. Hence, we can choose a ball with center $o$ and radius $t=\frac{1+a_0}{a_0}\mbox{diam}(\mathsf{A}(L)_s)$, which contains $H(u,\tau)\cap \mathsf{A}(L)$. In another words,
$$
H(u,\tau)\cap \mathsf{A}(L) \subseteq \mathsf{A}(L)_t.
$$
\par
For part 2, we can choose a number $\zeta>0$. By part 1, there exists a number $t>0$ such that, for any $u\in \omega$ and $H(u,\tau)\cap \mathsf{A}(L)_\zeta\neq \emptyset$, we have
$$
H(u,\tau)\cap \mathsf{A}(L)\subset \mathsf{A}(L)_t.
$$
\par
On the other hand, let $K\in \mathcal{C}_b(D)$ with $\mathsf{A}(K)=L$. By the definition of $K_\omega$, $K_\omega \cap \mathsf{A}(L)_\zeta\neq \emptyset$. Hence the supporting hyperplane $H(u,\tau)$ of $K_\omega$ with $u\in \omega$ satisfies $H(u,\tau)\cap \mathsf{A}(L)_\zeta\neq \emptyset$ and, therefore, $H(u,\tau)\cap \mathsf{A}(L)\subset \mathsf{A}(L)_t$ for some $t>0$, which gives $K_\omega^c\subseteq \mathsf{A}(L)_t$, hence $\mathsf{A}(L)_t\subseteq K_\omega$.
\epf

\bpf[Proof of Theorem \ref{M thm}]
First we define
$$
L:=\mathbb{H}^+_{n+1}\cap\left[\bigcap_{u\in \partial D}H^-(u,f(u))\right].
$$
It is clear that $L=\mathsf{A}(L)$. Given a nonzero finite Borel measure $\mu$ on $D$ with $\mu(D\backslash \omega)=0$, let $\mathfrak{C}(\omega)$ denote the set of positive continuous functions on $\omega$ with
$$f(u)
\begin{cases}
	>0 & u\neq -e_{n+1}\\
=0 & u=-e_{n+1}.
\end{cases}
$$
For any $f\in \mathfrak{C}(\omega)$, one has $K_{f}\in \mathcal{C}_c(D)$ satisfying
$$
K_{f}=L\cap \left(\bigcap_{u\in \omega} H^-(u,f(u))\right)
$$
Define a function $\Phi:\mathfrak{C}(\omega)\to (0,\infty)$ by
$$
\Phi(f):=V_n(K^c_{f})^{\frac{1}{n}}\int_\omega h_{K_f}(u)\,d\mu.
$$
The function $\Phi$ is continuous. We need to show that it attains a maximum on the set
$$
\mathfrak{L}:=\left\{h_M(u): M=L\cap \left(\bigcap_{u\in \omega} H^-(u,h_M(u))\right)\right\},
$$
Note that $o\in \mbox{int}(\omega)$. If we define
$$
C_{\omega}:=L\cap \left(\bigcap_{u\in \omega} H^-(u,0)\right),
$$
we always have $C_{\omega}\subseteq M$; hence, one has
$$
 V_n(L\backslash M)\leq V_n(L\backslash C_{\omega})=:\tau.
$$
\par
By Lemma \ref{lem5.1}, there is a number $t>0$ such that for any $h_M\in \mathfrak{L}$, one has $L\cap H(e_{n+1},t)\subseteq M$. This implies
$$
\Phi(h_M)\leq \tau^{n}\int_\omega h_{L_t}(u)\,d\mu(u)=:c,
$$
which is independent of $M$. Therefore,
$$
\sup\{\Phi(f): f\in \mathfrak{L}\}<\infty.
$$
\par
Let $\{K_i\}_{i\in \mathbb{N}}$ be a sequence with $h_{k_i}\in \mathfrak{L}$ such that
$$
\lim_{i\to \infty}\Phi(h_{K_i})=\sup\{\Phi(f): f\in \mathfrak{L}\}.
$$
For each $i$, we have $L\cap H(e_{n+1},t)\subset K_i$. By the Blaschke selection theorem, the bounded sequence $\{K_i\cap H^-(e_{n+1},t)\}_{i\in \mathbb{N}}$ of convex bodies has a subsequence converging to some convex body; hence, the same subsequence of $\{K_i\}_{i\in \mathbb{N}}$ converge to $K_0$, where $K_0\in \mathfrak{L}$.
By the continuity of $\Phi$, the function $\phi$ obtains the maximum at $h_{K_0}$. Note that $K_0\neq C_{\omega}$, otherwise, $\Phi(f)=0$ for all $f\in \mathfrak{C}(\omega)$.
\par
Therefore, let $f\in \mathfrak{C}(\omega)$. Then $h_{K_0}+\tau f\in\mathfrak{C}(\omega)$ for sufficiently small $|\tau|$; hence, the function
\be\label{equ5.20}
\tau\to V_n(K^c_{h_{K_0}+\tau f})^{\frac{1}{n}}\int_\omega h_{K_{h_{K_0}+\tau f}}(u)\,d\mu
\ee
obtains the maximum at $\tau=0$. Togeter with Lemma \ref{lem4.2}, we have the derivative of the function \ref{equ5.20} at $\tau=0$
$$
-\frac{1}{n}V_n(K_0^c)^{\frac{1-n}{n}}\int_\omega f\,dS(K,u)\int_\omega h_{K_0}\,d\mu+V_n(K_0^c)^{\frac{1}{n}}\int_\omega f\,d\mu(u)
$$
is equal to $0$. With
$$
c:=\frac{1}{nV_n(K_0^c)}\int_\omega h_{K_0}(u)\,d\mu(u),
$$
we have for any $f\in \mathfrak{C}(\omega)$,
$$
\int_\omega f\,d\mu(u)=c\int_\omega f\,dS(K_0,u)
$$
which gives
$$
d\mu=cdS(K_0,u).
$$
\epf
\section{The Minkowski problem for C-coconvex sets}
In the end, we show a solution to the Minkowski problem for $\sigma$-finite measure on $\Omega_C$.
\bpf[Proof of Theorem \ref{GM thm}]
Case 1, if $\mu(\omega)=0$ for any compact set $\omega\subseteq \Omega_C$, $K=C$ is the solution to the Minkowski problem for $\mu$, where $K$ is C-close set.
\par
Case 2, if $\mu(\Omega_C\backslash \omega)=0$ for some compact $\omega\subseteq \Omega_C$, the result will follow from \cite[Theorem 4]{Schneider18}.
\par
Case 3, we will use the idea in \cite[Section 12]{Schneider18}. If we choose a sequence $\{\omega_i\}_{i\in \mathbb{N}}$ of compact sets in $\Omega_C$ such that
$$
\nu(\omega_1)>0,\ \omega_i\subseteq \mbox{int}(\omega_{i+1}),\mbox{ and }\bigcup_{i}\omega_i=D.
$$
Since $\mu$ is a $\sigma$-finite measure, for each $j\in \mathbb{N}$, we define the measure $\mu_j$ by
$$
\mu_j(\eta):=\mu(\eta\cap \omega_i)\mbox{ for any Borel set } \eta\subseteq \Omega_C.
$$
Then $\mu_j$ is a nonzero, finite Borel measure with support concentrated on $\omega_j$. By \cite[Theorem 4]{Schneider18}, there exists a C-full set $K_j$ with $V_{n}(K_j^c)=1$ and 
$$
d \mu_j =c_jd S(K_j, \cdot)
$$
where
$$
c_j=\frac{1}{n}\int_{\omega_j}\bar{h}_{K_j}\,d\mu.
$$
We choose $t>0$ such that $V_n(C_t)\geq 1$, then
$$
K_j\cap C_t\neq \emptyset\mbox{ for all }j.
$$
Now taking same process in \cite[Section12]{Schneider18}, for an increasing sequence $\{t_k\}_{k\in \mathbb{N}}$ with $t_1>t$ and $t_k\to\infty$ as $k\to \infty$, there exists a subsequence such that
$$
K_{j_i}\cap C_{t_k}\to M_k\mbox{ as }i\to \infty,
$$
for each $k\in \mathbb{N}$. This gives
$$
K_i\cap C_{t_k}\to M_k \mbox{ as }i\to \infty\mbox{ for }k\in \mathbb{N}
$$
if changing $j_i$ by $i$, then $K_{j_i}=K_i$ and $\omega_{j_i}=\omega_i$. It is clear that for any $1\leq m<k$,
$$
M_m=M_k\cap C_{t_m}.
$$
Then define 
$$
M:=\cup_{k\in \mathbb{N}}M_k,
$$
then
$$
M\cap C_{t_k}=M_k.
$$
Hence, $M\subset C$ is a closed convex set.
\par
Since $V_n(K^c_i)=1$, we have
$$
V_n(M^c)=1.
$$
The last is to check that $c_i$ is convergent, where
$$
c_i=\frac{1}{n}\int_{\omega_i}\bar{h}_{K_j}\,d\mu.
$$
Note that the fact $V_n(M^c)=1$ gives there exists a $N\in \mathbb{N}$ such that $S(M,\omega_N)>0$; otherwise, $M=C$ gives $V_n(M^c)=0$. Since for each $i>N$, we have
$$
\mu(\omega_N)=\mu_i(\omega_{N})=c_iS(K_i,\omega_{N})
$$
and by \cite[Lemma 4.3]{YYZ}, for $i>N$,
$$
S(K_i,\omega_{N})\to S(M, \omega_{N})\mbox{ as }i\to \infty.
$$
hence, $S(K_i,\omega_{N})\geq \frac{1}{2}S(M,\omega_N)$ for sufficient large $i$. Then We set
$$
c:=\lim_{i\to \infty}\frac{\mu(\omega_{N})}{S(K_i,\omega_{N})},
$$
which implies
$$
c=\lim_{i\to \infty}\frac{1}{n}\int_{\omega_i}\bar{h}_{K_i}\,d\mu=\frac{1}{n}\int_{D}\bar{h}_M\,d\mu.
$$
In the end, for any compact set $\omega\subseteq D$, there exists $I\in \mathbb{N}$ such that $\omega\subseteq \omega_{I}$ with
$$
d\mu=c_id\bar{S}(K_i,\cdot) \mbox{ on }\omega
$$
for all $i>I$, which gives
$$
d\mu=cd\bar{S}(M,\cdot) \mbox{ on }\omega
$$
as $i\to \infty$.
\epf


   \vskip 2mm \noindent Ning Zhang, \ \ \ {\small \tt nzhang2@hust.edu.cn}\\
{ \em School of Mathematics and Statistics, Huazhong University of Science and Technology, 1037 Luoyu Road, Wuhan, Hubei 430074, China\\
Hausdorff Research Institute for Mathematics, Poppelsdorfer Allee 45
53115 Bonn, Germany} 
   
\end{document}